\input epsf

\documentclass[10pt]{amsart}

\usepackage{amscd, amsmath, amsthm}

\def\mylabel#1{\label{#1}   \rlap{\hskip1cm\leftline{#1}}   }
\def\mylabel#1{\label{#1}   \proplabeL{#1} \hskip-3pt  }

\def\mylabel#1{\label{#1}}

\def\pd{\partial}

\theoremstyle{definition}     

\numberwithin{equation}{section}

\newcommand{\vP}{{\mathaccent20 {\bf P}}}

\begin{document}


\title[Higher genus GW invariants of the Gr, Pf Calabi-Yau 
threefolds ]{ 
{Higher genus Gromov-Witten invariants of the Grassmannian, and 
the Pfaffian Calabi-Yau threefolds}
}

\author[S. Hosono]{ Shinobu Hosono }

\address{
Graduate School of Mathematical Sciences, 
University of Tokyo, Komaba Meguro-ku, 
Tokyo 153-8914, Japan
}
\email{hosono@ms.u-tokyo.ac.jp}

\author[Y. Konishi]{ Yukiko Konishi }

\address{
Graduate School of Mathematical Sciences, 
University of Tokyo, Komaba Meguro-ku, 
Tokyo 153-8914, Japan
}
\email{konishi@ms.u-tokyo.ac.jp}


\begin{abstract}
We solve Bershadsky-Cecotti-Ooguri-Vafa (BCOV) holomorphic anomaly 
equation to determine the higher genus Gromov-Witten invariants ($g \leq 5$) 
of the derived equivalent Calabi-Yau threefolds, which are of the 
appropriate codimensions in the Grassmannian $Gr(2,7)$ and 
the Pfaffian $Pf(7)$. 
\end{abstract}

\maketitle


\section{{\bf Introduction }}

Since the first successful application of the mirror symmetry to 
the Gromov-Witten theory of the quintic hypersurface in ${\bf P}^4$ 
\cite{CdOGP}, and its highly non-trivial generalization to higher genus 
$g \geq 1$ in \cite{BCOV1, BCOV2}, 
the mirror symmetry of Calabi-Yau manifolds 
has been attracting attentions in both mathematics and physics. Now, 
according to Kontsevich's homological mirror symmetry \cite{Ko}, we consider 
that two Calabi-Yau manifolds $X$ and $Y$ are mirror symmetric to each other 
when the derived category of the coherent sheaves on $X$, $D^b(Coh(X))$, 
is equivalent to the derived Fukaya category $DFuk(Y,\beta)$, and vice versa. 
In this homological viewpoint, it is clear that Calabi-Yau manifolds 
$X, X'$, which are derived equivalent 
$D^b(Coh(X)) \cong D^b(Coh(X'))$, are of considerable interest. 

For a smooth projective variety $X$, the projective varieties 
having equivalent derived category to $X$ are called Fourier-Mukai 
partners of $X$, and the set of their isomorphism classes is denoted by 
$FM(X)$. In dimension two, the set $FM(X)$ has been studied in detail 
in \cite{BM} and it has been shown that the number of Fourier-Mukai 
partners of a smooth minimal projective surface $X$ is finite, i.e. 
$|FM(X)|<\infty$. 
In particular, for a K3 surface $X$, a necessary and 
sufficient condition for a K3 surface $X'$ to be a partner 
of $X$ is known in terms of the Hodge isometry in the Mukai lattice 
\cite{Or}.  
Based on the result in \cite{Or}, a precise counting formula 
of the number of Fourier-Mukai partners has been given in 
\cite{Og, HLOY2}. 

In dimension three, however, since birational Calabi-Yau threefolds share the 
equivalent derived category \cite{Br}, the counting problem  should be 
considered under the birational equivalences. This contrasts with the 
fact in two dimensions that two birational K3 surfaces are 
biholomorphic to each other. Recently, an example of Calabi-Yau 
threefolds which share the equivalent derived category but seems 
to be non-birational has been constructed in \cite{BC, Ku} based on 
the earlier observation in \cite{Ro}. This example is of our  
interest in this article. 

In this paper, we apply the mirror symmetry to the derived 
equivalent Calabi-Yau threefolds $X$ and $X'$, that appeared 
in \cite{BC, Ro}, of appropriate codimensions 
in the Grassmannian $Gr(2,7)$ and the Pfaffian $Pf(7)$, respectively. 
In particular, we determine the higher genus Gromov-Witten invariants 
($g\leq 5$) integrating the holomorphic anomaly equation in \cite{BCOV2} 
recursively. The Gromov-Witten invariants at genus zero were determined 
earlier in \cite{BCFKvS} and \cite{Ro} following the method in \cite{CdOGP}. 
See \cite{Tj}\cite{BCFK}\cite{Ki} for mathematical proofs of the invariants. 
For the higher genus calculations, we solve the BCOV holomorphic anomaly 
equation \cite{BCOV1, BCOV2}. In particular we utilize the gap condition
at the conifold singularities, which has been found recently in \cite{HKQ}, 
with slight improvements in the estimate of the unknown parameters in 
the holomorphic ambiguities. 

Both the Calabi-Yau manifolds $X$ and $X'$ have Picard number $\rho=1$. 
Let $N_g^{X}(d)$, $N_g^{X'}(d)$ be the Gromov-Witten invariant of degree 
$d$ with respect to the respective generator $H$ of the Picard group. 
We denote by ${\tt F}_g(t)$ the generating functions of 
the Gromov-Witten invariants (, the so-called Gromov-Witten potential,) 
which have the following form for $g \geq 2$ in general,
\begin{equation}
{\tt F}_g(t)=
\frac{\chi}{2}\,(-1)^{g}\, \frac{|B_{2g} B_{2g-2}|}{2\,g\,(2g-2)\,(2g-2)!} 
+ \sum_{d > 0} N_g(d) \, q^d \;\;, \;\;(q=e^{2\pi i t})
\mylabel{eqn:Fg-def}
\end{equation}
where $\chi$ is the Euler number of a Calabi-Yau 
manifold and $B_g$ is the Bernoulli number. The constant term above 
represents the Gromov-Witten invariant $N_g(0)$ of degree zero, 
and it represents 
the contribution from the constant maps \cite{BCOV2, FP}. We determine 
the potential ${\tt F}_g(t)$ for $g \leq 5$ for $X$ and $X'$. 
Also to see some implications of our results to the enumerative problem 
of holomorphic curves and/or the moduli problem related to 
Donaldson-Thomas invariants \cite{T}, we list the so-called 
Gopakumar-Vafa `invariants' $n_g(d)$ \cite{GV} which are 
determined from $N_g(d)$ by 
\begin{equation}
\sum_{g \geq 0} \lambda^{2g-2} {\tt F}_g(t) = 
\sum_{g\geq 0} \sum_{k \geq 1,d \geq 0} 
n_g(d) \, \frac{1}{k}\,\Big(2 \,\sin \frac{k\lambda}{2} \Big)^{2g-2}\,q^{kd} \;\;.
\mylabel{eqn:GV-inv}
\end{equation}

\vskip0.3cm
The organization of this paper is as follows: In section 2, we summarize the 
constructions of the Grassmannian and the Pfaffian Calabi-Yau threefolds, and 
their mirror orbifolds given in \cite{Ro}.  After introducing the 
Picard-Fuchs differential equation of the period integrals, we 
determine the $g=0,1$ Gromov-Witten prepotentials. We will also make a 
comment on a similarity of the Picard-Fuchs differential equation to 
the corresponding differential equation studied for a K3 surface with 
a non-trivial Fourier-Mukai partner in \cite{HLOY2}. 
In section 3, we briefly introduce the BCOV holomorphic anomaly equation for 
$g\geq 2$ and its solutions given in \cite{BCOV2}. According to 
\cite{BCOV1,BCOV2} we define the topological limit of the solutions. 
We also summarize the recent results found in \cite{YY} about some 
polynomiality of the solutions. Then, we introduce a `gap condition', 
which has been found recently in \cite{HKQ}, to fix the holomorphic 
ambiguities contained in the solutions. 
In section 4, we present our calculations in some details to determine 
the Gromov-Witten potentials. We determine the potentials up to $g=5$ and 
list the resultant Gopakumar-Vafa invariants in Tables 1 and 2. The conclusion 
and discussions are given in section 5.

\vskip0.3cm
\noindent
{\bf Acknowledgments:}
The authors would like to thank M.-H. Saito for letting them know the 
reference \cite{BH}. They are grateful to  I. Ciocan-Fontanine and B. Kim 
for explaining the results in \cite{BCFK}\cite{Ki}.  
S. H. would like to thank C. Doran for drawing his 
attention to the differential equation of the Grassmannian 
and the Pfaffian Calabi-Yau manifolds. 
The work by S. H. is supported in part by Grant-in Aid Scientific Research 
(C 18540014). The work by Y. K. is supported in part by JSPS Research 
fellowships for Young Scientists.

\vskip1cm


\section{{\bf Gromov-Witten invariants at $g=0$ and $g=1$} }

In this section, we briefly summarize the constructions of 
Calabi-Yau manifolds, $X$ and $X'$, and their orbifold mirror construction 
following \cite{Ro}. We summarize the genus zero and one Gromov-Witten 
invariants using the solutions of the Picard-Fuchs differential equation 
of the mirror family.

\vskip0.5cm
\noindent
{\bf (2-1) The Grassmannian and the Pfaffian Calabi-Yau threefolds.}
Let us first summarize the construction of the Grassmannian Calabi-Yau 
threefold and its topological invariants. 
Let $Gr(2,7)$ be the Grassmannian of the 2-planes in ${\bf C}^7$, and 
$Q$ be the universal quotient bundle. The line bundle $\wedge^5 Q$ 
determines the Pl\"ucker embedding $i: Gr(2,7) \hookrightarrow {\bf P}^{20}$, 
hence $\sigma_1=c_1(Q)$ represents the class of a hyperplane section. Then 
$\int_{Gr(2,7)}\sigma_1^{10}=42$ gives the degree of the Grassmannian in 
the projective space. We denote by $Gr(2,7)_{1^7}$ the complete 
intersection of $Gr(2,7)$ with seven hyperplanes in ${\bf P}^{20}$. Then 
$X=Gr(2,7)_{1^7}$ defines a Calabi-Yau threefold since 
$c_1(Gr(2,7))=7 \sigma_1$. In fact, the Chern class of $Gr(2,7)$ is 
expressed by  
\[
c(Gr(2,7))=1+7 c_1(Q)+(25 c_1(Q)^2-3 c_2(Q))+14(4 c_1(Q)^3-c_1(Q)c_2(Q))+ 
\cdots ,
\]
see e.g. \cite{BH}, and for the complete intersection, we have
\[
c(X)=\frac{c(Gr(2,7))}{(1+c_1(Q))^7} = 
1+(4c_1(Q)^2-3c_2(Q)) -7(c_1(Q)^3-c_1(Q)c_2(Q)) \;\;.
\]
Using $\int_{Gr(2,7)}\sigma^{10}_1=42$, 
$\int_{Gr(2,7)}\sigma_1^8 \sigma_2 =
28$ $(\sigma_2=c_2(Q))$, and representing by $H$ the hyperplane $\sigma_1$ 
on $X$, we have the following topological invariants
\begin{equation}
\chi(X)= -98 \;, \;\; c_2(X)\cdot H = 84 \;, \;\; H^3=42
\mylabel{eqn:topGr}
\end{equation}
Also we see $h^{1,1}(X)=1$ by Lefschetz hyperplane theorem, which 
implies $h^{2,1}(X)=50$.

\vskip0.3cm

The construction of the second Calabi-Yau manifold $X'$ is more involved and 
utilizes the Pfaffian variety in the projective space ${\bf P}^{20}$. 
Let $\mathcal S$ be a $7 \times 7$ skew symmetric matrix 
${\mathcal S}=(s_{ij})$ 
with $[s_{ij}] \in {\bf P}^{20}$. The rank of $\mathcal S$ is less than 
or equal to 6, and in particular, the rank 4 locus 
(${\rm rk} \mathcal S \leq 4$) determines a
codimension three variety in ${\bf P}^{20}$, the Pfaffian variety. 
Explicitly, this variety is determined by 
the ideal generated by the square roots of the diagonal minors of $\mathcal S$,
$p_0(\mathcal S), \cdots, p_6(\mathcal S)$. Restricting this variety to 
a generic projective 
space ${\bf P}^6 \subset {\bf P}^{20}$, i.e. specializing the parameters 
$[s_{ij}]$ to lie on a generic ${\bf P}^6$, we have an exact sequence
\[
0 \rightarrow 
{\mathcal O}_{{\bf P}^6}(-7) 
\xrightarrow[]{\,^t p(\mathcal S)} 
{\mathcal O}_{{\bf P}^6}(-4)^{\oplus 7}  
\xrightarrow[]{ {\mathcal S} } 
{\mathcal O}_{{\bf P}^6}(-3)^{\oplus 7}  
\xrightarrow[]{ p(\mathcal S) } 
{\mathcal O}_{{\bf P}^6}  
\rightarrow 
{\mathcal O}_{X'} 
\rightarrow 0 \;,
\]
where we set $p(\mathcal S)=( \, (-1)^{i+1}p_i(\mathcal S) \,)_{i=0,..,6}$ 
to be a row vector and use $p(\mathcal S)\,\mathcal S = 0$, 
since $(-1)^{i+j} p_i(\mathcal S)p_j(\mathcal S)$ represents the $ij$-minor 
of $\mathcal S$ and ${\rm det}(\mathcal S)=0$. From this exact sequence, 
we see that the canonical sheaf of $X'$ is trivial, 
$\omega_{X'} \cong 
{\mathcal Ext}^3( {\mathcal O}_{X'}, \omega_{{\bf P}^6}) 
\cong {\mathcal O}_{X'}$ and therefore $X'$ is a Calabi-Yau threefold. 
 The degree of $X'$ in ${\bf P}^6$ is 14 and $H^3=14$ for the hyperplane 
section $H$. Other topological invariants are determined by the general 
formulas $c_2(X')H=84-2d, c_3(X')=-d^2+49d-588$ $(d=H^3)$ valid for 
codimension 3 smooth Calabi-Yau varieties, see \cite{To} for example. 
Thus we have the following topological invariants
\begin{equation}
\chi(X')= -98 \;, \;\; c_2(X')\cdot H = 56 \;, \;\; H^3=14 \;\;.
\mylabel{eqn:topPf}
\end{equation}
We also have the Hodge numbers $h^{1,1}(X')=1$ and $h^{2,1}(X')=50$.

As noted in the reference \cite{Ro}, the construction of $X$ and $X'$ are 
dual in the following sense,
\[
X = Gr(2,7) \cap {\bf P}^{13} \subset {\bf P}^{20}  \;\;,\;\;
X' = Pf(7)  \cap \vP^6 \subset \vP^{20},
\]
where $\vP^{20}$ is the dual projective space to ${\bf P}^{20}$ and 
$\vP^6$ is the annihilator of ${\bf P}^{13}$ under the dual pairing. 
This duality has been utilized 
to prove the derived equivalence $D^b(Coh(X)) \cong D^b(Coh(X'))$ 
\cite{BC,Ku}.

\vskip1cm
\noindent
{\bf (2-2) The mirror manifolds and the Picard-Fuchs differential equations.}  
In a similar way to the orbifold construction of the quintic hypersurface 
in ${\bf P}^4$ \cite{GP, CdOGP}, the mirror manifolds $Y$ and $Y'$ 
have been constructed, respectively, for $X$ and $X'$ in \cite{Ro}.  
Following \cite{Ro}, we introduce the mirror family 
$\{ Y'_x \}_{x \in {\bf P}^1}$ and the Picard-Fuchs differential equation 
for the period integral of a holomorphic three form. 

\vskip0.3cm
Consider the skew symmetric matrices 
\[
E_k(y)=\sum_{i+j=k} y_{i-j} E_{ij} \;\;\;(k=0,1,\cdots,6; \; y_i+y_{-i}=0)
\]
parametrized by $[y_1,y_2,y_3] \in {\bf P}^2$, where the index of $y_i$ is 
understood modulo 7, and $E_{ij} (0 \leq i,j \leq 6)$ are the matrix units. 
We define 
\[
\vP^6_{[y_1,y_2,y_3]}=\text{the projective span of } 
\{ E_0(y), \cdots, E_6(y) \} 
\subset \vP^{20} \;,
\]
and consider the group $G=\langle \tau,\sigma \rangle$ acting 
on $\vP^6_{[y_1,y_2,y_3]}$ by 
\[
\tau: E_k(y) \mapsto e^{2\pi i k/7} E_k(y) \;\;,\;\;
\sigma: E_k(y) \mapsto E_{k+2}(y) \;\;. 
\]
Then $X'_{[y_1,y_2,y_3]}=Pf(7) \cap \vP^6_{[y_1,y_2,y_3]}$ is a special 
family of $X'$, and its general member has 49 double points at 
the orbit $G \cdot [y_0,y_1,\cdots,y_6]$. 
When we further specialize $X'_{[y_1,y_2,y_3]}$ 
to a ${\bf P}^1$ family $X'_{[y_1,y_2,0]}$, we have additionally 7 
double points at the fixed points of $\tau$. These double points 
arise from the process collapsing $S^3$ to points, and may be blown 
up to ${\bf P}^1$ without affecting the Calabi-Yau condition. 
Blowing up these $49+7$ double points in total, we have 
$-98+(49+7)\times 2 =14$ for the Euler number of the resolved 
space $\tilde X'_{[y_1,y_2,0]}$. Now, consider the quotient 
$\tilde X'_{[y_1,y_2,0]}/ \langle \tau \rangle$. This quotient has 
singularities which come from the $7 \times 2$ fixed points under 
the action of $\tau$. These singularities can be resolved under the Calabi-Yau 
condition, and for the Euler number we have 
\[
\chi(\widehat{\tilde X'_{[y_1,y_2,0]}/\langle \tau \rangle} )=
\frac{1}{|\langle \tau \rangle |} 
\sum_{g,h \in \langle \tau \rangle} 
\chi(\tilde X'_{[y_1,y_2,0]}|_{{g,h}}) = 98 \;\,
\]
where $\tilde X'_{[y_1,y_2,0]}|_{g,h}$ 
represents the fixed points under $g$ and $h$ (, i.e.  
the $7 \times 2$ points for $(g,h)\not=(e,e)$). 
The Hodge numbers are determined in \cite{Ro} by looking 
the blow-ups more closely. The results are $h^{1,1}=50, h^{2,1}=1$, 
justifying the claim that 
$Y'_y:=\widehat{\tilde X'_{[y_1,y_2,0]}/\langle \tau \rangle }$ is the mirror 
family of the Calabi-Yau variety $X'$. 

\vskip0.5cm

For the concrete description of $Y'_y$, we write in Appendix (A-1) the 
diagonal Pfaffians $p_k(\mathcal S)$ of the skew symmetric matrix 
$\mathcal S(y)=\mathcal S(y,[u])$ for 
the spacial family $X'_{[1,y,0]}$ with 
$[u]=[u_0,u_1,\cdots,u_6] \in {\bf P}^6$. From the explicit form of the 
generators $p_k(\mathcal S)$, we see that 
$Y'_{e^{2\pi i/7}y} \cong Y'_y$ 
and hence $x=y^7$ parametrize the genuine mirror family. 
Then in terms of $p_i(\mathcal S)$, the 
holomorphic three form of the family $Y'_x$ may be given by 
\begin{equation}
\Omega(x)={\rm Res} 
\frac{(-1)^{\epsilon} 
Pf({\mathcal S}_{i_3i_4i_5i_6}) d\mu }{p_{i_0}p_{i_1}p_{i_2}} \;\;,
\mylabel{eqn:Omega-Res}
\end{equation}
where $d\mu=du_0du_1\cdots du_6$ and 
${\mathcal S}_{i_3i_4i_5i_6}$ is the $4 \times 4$ `diagonal' 
sub-matrix of $\mathcal S$ 
specified by the index set $\{ i_3i_4i_5i_6 \}$, and 
$\epsilon$ represents the parity of the order $i_0i_1\cdots i_6$.   
Evaluating the period integral over a 
torus cycle as a power series in $x$, the Picard-Fuchs 
differential operator ${\mathcal D}_x$ has been 
determined in \cite{Ro}, 
\begin{equation}
\begin{aligned}
{\mathcal D}_x=&9 \,\theta_x^4 - 
3\,x(15+102\,\theta_x+272\,\theta_x^2+340 \,\theta_x^3+173\,\theta_x^4) \\
&-2\,x^2(1083+4773\,\theta_x+7597\,\theta_x^2+5032\,\theta_x^3+1129\,\theta_x^4)\\ 
&+2\,x^3(6+675\,\theta+2353\,\theta_x^2+2628\,\theta_x^3+843\,\theta_x^4) \\
&-x^4(26+174\,\theta_x+478\,\theta_x^2+608\,\theta_x^3+295\,\theta_x^4)+
x^5(\theta_x+1)^4 \;,
\end{aligned}
\mylabel{eqn:GrPF}
\end{equation}
where we define $\theta_x=x\frac{d \;}{dx}$. Using this differential operator, 
and normalizing the holomorphic three form suitably, 
we can determine the Yukawa coupling to be 
\begin{equation}
C_{xxx}:=\int_{Y'_x} \Omega(x) \wedge \left(\frac{d \;}{dx}\right)^3 \Omega(x) 
= \frac{42-14x}{x^3(1-57x-289x^2+x^3)} \;.
\mylabel{eqn:Cxxx}
\end{equation}

\vskip0.5cm

A similar orbifold construction works for the Grassmannian Calabi-Yau 
variety $X$ by taking the dual projective space 
${\bf P}_{[1,y,0]}^{13}$ to $\vP^6_{[1,y,0]}$. 
Then the mirror family $\{ Y_y \}$ 
is given by a resolution of a suitable orbifold of 
$Gr(2,7)\cap {\bf P}^{13}_{[1,y,0]}$. The remarkable observation made 
in \cite{Ro} is that we obtain the same Picard-Fuchs differential operator 
as above, which has the property of the maximally degeneration \cite{Mo} 
at both $x=0$ and $x=\infty$. 
This indicates that the two Calabi-Yau varieties $X$ and $X'$ 
share the same the mirror family $\{Y_y\} = \{Y'_y\}$, 
and that the complexified 
K\"ahler moduli of the two Calabi-Yau varieties $X$ and $X'$ are 
unified in one complex structure moduli of the mirror family 
(Conjecture 10 in \cite{Ro}). The structure of the singularities of the 
Picard-Fuchs equation (\ref{eqn:GrPF}) (cf. \cite{DM, ES})  
may be summarized in the following 
Riemann's P scheme listing the indices $\rho_k$ at each regular singular point;
\begin{equation}
\left\{
\begin{array}{c|cccccc}
x     & 0 & \alpha_1 & \alpha_2 & \alpha_3 & 3 & \infty \\
\hline
\rho_1 & 0 &    0    &   0    &   0    &  0  &  1      \\
\rho_2 & 0 &    1    &   1    &   1    &  1  &  1      \\
\rho_3 & 0 &    1    &   1    &   1    &  3  &  1      \\
\rho_4 & 0 &    2    &   2    &   2    &  4  &  1      \\
\end{array}
\right\} \;\;,
\mylabel{eqn:P-indices}
\end{equation}
where $\alpha_k$ are the roots of the 'discriminant' $1-57x-289x^2+x^3=0$, 
for which $Y'_x$ has double points. 

Making the instanton expansions at each 
degeneration point, we find that the expansion about $x=0$ 
corresponds to the K\"ahler moduli of the Grassmannian Calabi-Yau $X$, 
and the expansion about $x=\infty$ to that of the Pfaffian Calabi-Yau $X'$. 
Our main objective in this paper is to extend the instanton 
calculations to higher genera. 

\vskip1cm
\noindent
{\bf (2-3) A digression to K3 surfaces.} 
It is clear that the property of the Picard-Fuchs differential operator 
${\mathcal D}_x$ 
is closely related to the equivalence 
$D^b(Coh(X))$ $\cong$ $D^b(Coh(X'))$. 
Here we remark that essentially the same property 
may be observed in the case of K3 surfaces. 

Let us recall that the set of Fourier-Mukai partners for a smooth projective 
variety $X$ is defined by 
\[
FM(X) = \{ Y \;|\; D^b(Coh(Y)) \cong D^b(Coh(X)) \;\}/ \sim \;\;,
\]
where $\sim$ represents the isomorphisms. When $X$ is a K3 surface, 
one may expect that the cardinality of $|FM(X)|$ is finite since 
birational K3 surfaces are biholomorphic. In fact, it is known that 
the number of Fourier-Mukai partners is finite\cite{BM, Or}. 
In particular, for a K3 surface $X$ of degree $2n$ and the Picard 
number $\rho(X)=1$, it is found \cite{Og} that the number of 
the Fourier-Mukai partners has a simple form,
\[
|FM(X)|=2^{p(n)-1} \;\;,
\]
where $p(n)$ is the number of the prime factors ($p(1):=1$). The first 
nontrivial case arises from $p(6)=2$, i.e., we have $FM(X)=\{ \, X, X' \,\}$ 
for a K3 surface $X$ of degree 12. According to \cite{Mu}, the partner $X'$ 
may be identified with a moduli space of the rank 2 stable sheaves with 
$c_1(\mathcal E)=H, \chi(\mathcal E)=2+3$. Also  
explicit constructions of the K3 surfaces of degree 12 and the mirror 
K3 surfaces are known in detail, see \cite{HLOY1} and references therein. 
There the modular group $\Gamma(6)_{0+}$ appears as the monodromy group 
of the Picard-Fuchs differential equation of the mirror (one-parameter) 
family. It is found in \cite{HLOY1} that one of the generators of 
the group $\Gamma(6)_{0+}$ does not correspond to any element in  
$Auteq(D^b(Coh(X)))$ under the mirror symmetry, and argued that 
this generator represents the Fourier-Mukai transform $\Phi^{\mathcal P}$ 
defined by the Poincar\'e bundle ${\mathcal P}$ over $X \times X'$. 
The rest of the generators defines the index two subgroup 
$\Gamma(6)_{0+6}$ of $\Gamma(6)_+$.  

Changing the monodromy group to the smaller group $\Gamma(6)_{0+6}$ 
doubles the moduli space of the mirror family (or the fundamental 
domain in the upper half plane). 
This doubled mirror family may be found in the table of \cite{LY}, and it has 
the Picard-Fuchs differential operator, 
\[
{\mathcal D}=\theta_x^3-x(2\,\theta_x+1)(17\,\theta_x^2+17\,\theta_x+5) 
+x^2(\theta_x+1)^3 \;\;. 
\]
This differential operator shows exactly the same property as 
(\ref{eqn:GrPF}), i.e.,  it has two maximal degeneration points 
at $x=0$ and $x=\infty$. One may pursue the similarity further in 
that the Fourier-Mukai partner $X'$ has an explicit 
construction using the orthogonal Grassmannian \cite{Mu}. 
Here, a naive construction of the Grassmannian K3 surface 
$X=Gr(2,6)_{1^6}$, however, does not give ${\rm deg}(X)=12$ 
but $14$, and hence $|FM(X)|=1$.

\vskip1cm
\noindent
{\bf (2-4) $g=0$ and $g=1$ Gromov-Witten invariants.}  
We summarize the calculations of the genus zero and 
one Gromov-Witten invariants of the Grassmannian and the Pfaffian 
Calabi-Yau varieties $X,X'$.

\vskip0.5cm
\noindent
{\bf (2-4.\,a)} 
Let us first introduce the so-called mirror map \cite{CdOGP}.  
We will denote henceforth the local coordinate $z=\frac{1}{x}$ to analyze the 
local solutions of the Picard-Fuchs equation (\ref{eqn:GrPF}) about 
$x=\infty$. At each degeneration point, we have one regular series solution 
with other solutions having (higher) logarithmic singularities. 
We normalize the regular solution and choose the following 
linear-logarithmic solution 
\begin{equation}
\begin{cases}
w_0(x)= 1 + 5\,x + 109\,x^2 +
3317\,x^3 + 121501\,x^4 + \cdots \;, \\
w_1(x)=\log(x) w_0(x)+14\,x + 357\,x^2 + \frac{35105}{3}\,x^3 + 
\frac{2669975}{6}\,x^4 + \cdots \;\;. 
\end{cases}
\mylabel{eqn:w0(x)}
\end{equation}
The choice of the linear-logarithmic solution $w_1(x)$ is up to the addition 
of arbitrary multiple of $w_0(x)$.  Here we fix this ambiguity so that 
the complexified K\"ahler moduli $2 \pi i t=\frac{w_1(x)}{w_0(x)}$ has a 
`nice' form of the $q$-expansion,
\[
\frac{1}{x(q)}=\frac{1}{q}+ 14 + 189\,q + 2534\,q^2 + 42826\,q^3 + 869162\,q^4 + 
\cdots \;\;,
\]
where $q:=e^{2\pi i t}$. In a similar way, we fix the regular solution 
$\tilde w_0(z)$ and the linear-logarithmic solution $\tilde w_1(z)$ at $z=0$,
\begin{equation}
\begin{cases}
\tilde w_0(z)=z + 17\,z^2 + 1549\,z^3 + 215585\,z^4 + 36505501\,z^5 + 
 \cdots\;, \\
\tilde w_1(z)= \log(z) \tilde w_0(z) + 
70\,z^2 + 7413\,z^3 + \frac{3268573\,z^4}{3} +
\frac{1138372375\,z^5}{6} + \cdot\cdot .
\end{cases}
\mylabel{eqn:w0(z)}
\end{equation}
By defining $2 \pi i \tilde t = \frac{\tilde w_1(x)}{\tilde w_0(x)}$, 
$\tilde q=e^{2\pi i \tilde t}$, we have 
\[
\frac{1}{z(\tilde q)}= \frac{1}{\tilde q} + 
   70 + 3773\,\tilde q + 232750\,\tilde q^2 + 18421802\,\tilde q^3 + 
   1781859058\,\tilde q^4 + \cdots \;.
\]
The expansions $x=x(q)$ and $z=z(\tilde q)$ are called mirror maps 
at the respective degeneration points, $x=0$ and $z=0$. 

\vskip1cm
\noindent
{\bf (2-4.\,b)}
Now, by the formula in \cite{CdOGP}, we determine the quantum corrected 
Yukawa coupling $K_{ttt}(t)$ at $x=0$ by 
\[
\Big(\frac{1}{w_0(x)}\Big)^2 C_{xxx} \Big(\frac{d x}{dt}\Big)^3 
=  42 + 196\,q + 9996\,q^2 + 344176\,q^3 + 12685708\,q^4 +  \cdots \;.
\]
For the expansion at $z=0$, we transform the Yukawa coupling (\ref{eqn:Cxxx}) 
by
\begin{equation}
C_{zzz}(z)=C_{xxx}(x)\Big(\frac{dx}{dz}\Big)^3 = 
\frac{14-42\,z}{z(1 - 289\,z - 57\,z^2 + z^3)} \;\;.
\mylabel{eqn:Czzz}
\end{equation}
Then the quantum Yukawa coupling $K_{\tilde t \tilde t \tilde t}(\tilde t)$ 
at $z=0$ is given by 
\[
\Big(\frac{1}{\tilde w_0(x)}\Big)^2 C_{zzz} \Big(\frac{d z}{dt}\Big)^3 
=   14 + 588\,q + 97412\,q^2 + 15765456\,q^3 + 2647082116\,q^4 + 
\cdots .
\]
These Yukawa couplings are related to the Gromov-Witten potentials by
\[
K_{ttt}(t)=\Big(q \frac{d \;}{d q} \Big)^3 {\tt F}_0(t)  \;\;,\;\;
K_{\tilde t \tilde t \tilde t}(\tilde t)=
\Big(\tilde q \frac{d \;}{d \tilde q} \Big)^3 \tilde {\tt F}_0(\tilde t) \;.
\]
Comparing the topological data given in (\ref{eqn:topGr}) and 
(\ref{eqn:topPf}), the degenerations at $x=0$ and $z=0$ have been identified  
in \cite{Ro}, respectively, with the Grassmannian Calabi-Yau $X$ and 
the Pfaffian Calabi-Yau $X'$, i.e., 
\[
{\tt F}_0(t) = {\tt F}_0^X(t) \;\;,\;\;
\tilde{\tt F}_0(\tilde t) = {\tt F}_0^{X'}(\tilde t) \;.
\]
We observe in (\ref{eqn:Czzz}) that the numerator 
of the Yukawa coupling, i.e. $42-14\,x$, explains the 
difference of the leading  term between the $q$- and the 
$\tilde q$-expansions. This simple observation should be contrasted to the 
similar calculations done for the `topology changes' (, i.e. flops)\cite{AGM}.

\vskip1cm
\noindent
{\bf (2-4.\,c)} 
For the genus one invariants, we apply the BCOV formula\cite{BCOV1} of the 
holomorphic potential $F^{(1)}(x)$ to our case, 
\begin{equation}
F^{(1)}(x)=\frac{1}{2}
\log\Big\{ \Big(\frac{f_1(x)}{w_0(x)}\Big)^{3+h^{1,1}-\frac{\chi}{12}}
\Big( \frac{d x}{dt} \Big)  \, dis(x)^{-\frac{1}{6}} \, x^{-1-\frac{c2.H}{12}} 
\Big\} \;,
\mylabel{eqn:F1}
\end{equation}
where $dis(x)=1 - 57\,x - 289\,x^2 + x^3$ and $f_1(x)$ is some holomorphic 
function which we fix to $f_1(x)=1$ by requiring the regularity of 
$F^{(1)}(x)$ at $x=0,\infty,3$. 
Exactly the same form as $F^{(1)}(x)$ applies to $\tilde F^{(1)}(z)$ with 
$\tilde w_0(z)$, $\widetilde{dis}(z)=1 - 289\,z - 57\,z^2 + z^3$, 
$\tilde f_1(z)=z$ and the data (\ref{eqn:topPf}). The 
holomorphic function $\tilde f_1(z)$ guarantees  the regularity of 
$\tilde F^{(1)}(z)$ at $z=0$. 
Using the topological data (\ref{eqn:topGr}),(\ref{eqn:topPf}) and the mirror 
maps $x=x(q)$ and $z=z(\tilde q)$, we obtain the genus one Gromov-Witten 
potentials,
\[
{\tt F}_1^X(t)=F^{(1)}(x(q)) \;\;,\;\;
{\tt F}_1^{X'}(\tilde t)=\tilde F^{(1)}(z(\tilde q)) \;\;.
\]

Here we remark that, except that one has to replace $w_0(x)$ 
with $\tilde w_0(z)$ by hand, one can verify the equality 
\[
F^{(1)}(x) = \tilde F^{(1)}(z)  \;\;,
\]
with $x=\frac{1}{z}$. This relation holds because, 
by taking the topological limits,  
the BCOV formulas (\ref{eqn:F1}) and $\tilde F^{(1)}(z)$ follow 
from the `Quillen's norm' function 
\begin{equation}
{\mathcal F}^{(1)}(x,\bar x)=\frac{1}{2} \log 
\Big\{ e^{(3+h^{11}-\frac{\chi}{12})\,K}\,G^{x\bar x} 
\Bigl|dis(x)^{-\frac{1}{6}} \, x^{-1-\frac{c2.H}{12}}\Bigr|^2 \Big\} \;,
\mylabel{eqn:mathcal-F1}
\end{equation}
of a certain holomorphic bundle over the moduli space \cite{BCOV2}, 
see also \cite{BT, FZY}. 
We will define the topological limits in (\ref{eqn:top-limit}) and come to 
this point in the next section, see subsection (3-5).

\newpage
\section{{\bf BCOV holomorphic anomaly equation}}

Here we introduce the BCOV holomorphic anomaly equation and its 
topological limits at the degeneration points in the moduli space. 

\vskip0.5cm
\noindent
{\bf (3-1) The special K\"ahler geometry. }
The mirror family $\{ Y_x \}_{x \in {\bf P}^1}$ defines the so-called 
special K\"ahler geometry on each neighborhood $B_0$ of $x_0 (\not= 0, 
\alpha_1,\alpha_2,\alpha_3,\infty)$ on the moduli space ${\bf P}^1$. 
Let us denote ${\mathcal M}^{cpl}={\mathbf P}^1\setminus
\{0,\alpha_1,\alpha_2,\alpha_3,\infty\}$. 
To describe the geometry on ${\mathcal M}^{cpl}$, let 
$\Omega(x)=\Omega(Y_x)$ $(x \in B_0)$ be the holomorphic 
three form (\ref{eqn:Omega-Res}), normalized by (\ref{eqn:Cxxx}). Consider 
the middle cohomology $H^3_{x_0}=H^3(Y_{x_0},{\bf Z})$, and define 
the period domain,
\[
D=\{ \; \omega \in {\bf P}(H^3_{x_0}\otimes {\bf C}) \;|\; 
(\omega, \omega)=0 \;,\; (\omega, \bar\omega) >0 \} \;,
\]
where $(\omega,\omega'):=i \int_{Y_{x_0}} \omega \wedge \omega'$. 
Making an identification $H^3(Y_x,{\bf Z}) \cong H^3(Y_{x_0},{\bf Z})$ 
for $x \in B_0$, the choice of the holomorphic three form $\Omega(x)$ 
determines the period map ${\mathcal P}_0: B_0 \rightarrow D$.  
Let ${\mathcal U}$ be the restriction of the tautological line bundle 
of ${\bf P}(H^3_{x_0}\otimes {\bf C})$ to $D$. Then we have a holomorphic 
line bundle ${\mathcal L}={\mathcal P}^*_0{\mathcal U}$ over $B_0$. 
Globalizing 
this local construction, we obtain a holomorphic line bundle ${\mathcal L}$ 
over a covering space $\tilde {\mathcal M}^{cpl}$ with its covering 
group (`modular group') $\Gamma \subset Sp(4,{\mathbf Z})$. 

The special K\"ahler geometry on $B_0$ is defined by the Weil-Peterson 
metric $G_{x\bar x}=\pd_x \pd_{\bar x} K(x,\bar x)$ with 
the K\"ahler potential $K(x,\bar x)= - \log (\Omega(x),\overline{\Omega(x)})$. 
Since $K(x,\bar x)$ is monodromy invariant, we see that this local 
geometry naturally glues together on  
${\mathcal M}^{cpl}$. Consider the metric connection given by 
$\Gamma_{xx}^{\;\;x}=G^{x\bar x} \pd_x G_{x\bar x}$ and 
$\Gamma_{\bar x \bar x}^{\;\;\bar x}=G^{\bar x x} 
\pd_{\bar x} G_{\bar x x}$. 
This connection defines the covariant derivative on 
the sections of the tangent bundle $T{\mathcal M}^{cpl} \otimes {\bf C} 
= T'{\mathcal M}^{cpl} \oplus T''{\mathcal M}^{cpl}$. Then we may 
write the so-called special K\"ahler geometry relation, 
\begin{equation}
\pd_{\bar x} \Gamma_{xx}^{\;\;x} = 2 G_{x\bar x} - 
C_{xxx} C_{\bar x \bar x \bar x}e^{2K} G^{x \bar x} G^{x \bar x} \;,
\mylabel{eqn:sp-geom-rel}
\end{equation}
where $K$ is the K\"ahler potential and $C_{xxx}$ is the 
Yukawa coupling (\ref{eqn:Cxxx}). It is known that this relation 
follows from a certain local system over ${\mathcal M}^{cpl}$ associated to 
$H^3(Y_x,{\mathbf Z})$\cite{St1}.

Now let us introduce `K\"ahler connection' by $K_x=\pd_x K$ and 
$K_{\bar x}=\pd_{\bar x} K$. We see that this connection 
defines the covariant derivative on the sections of ${\mathcal L}$, more 
precisely `$\Gamma$-modular forms of weight one', and also its complex 
conjugate $\bar{\mathcal L}$, and the tensor products thereof.  
We have $D_x \xi = \pd_x \xi + n K_x \xi + m K_{\bar x} \xi $ 
for a section $\xi \in \Gamma({\mathcal L}^{n} \otimes \bar{\mathcal L}^{m})$, 
again more precisely, for a `$\Gamma$-modular form' $\xi$ of weight $(n,m)$. 
Thus for a holomorphic 
tangent vector $\xi^x$ taking a value in $\bar{\mathcal L}$, for example, 
we have 
$D_{x} \xi^x = ( \pd_x + \Gamma_{xx}^{\;\; x} ) \xi^x$ and 
$D_{\bar x} \xi^x = ( \pd_{\bar x} + K_{\bar x} ) \xi^x$.

\vskip1cm \noindent
{\bf (3-2) BCOV anomaly equation and the general solutions ${\mathcal F}^{(g)}$. }
Using the special K\"ahler geometry and also the Griffiths transversality 
for the period map, we can show that there exist potential functions 
which express the Yukawa coupling (\ref{eqn:Cxxx}) and its complex conjugate 
by 
\[
C_{xxx}=D_xD_xD_x {\mathcal F}^{(0)}(x,\bar x)  \;\;, \;\; 
C_{\bar x \bar x \bar x}=D_{\bar x}D_{\bar x}D_{\bar x} 
\bar {\mathcal F}^{(0)}(x,\bar x)  \;\;, 
\]
where ${\mathcal F}^{(0)}(x,\bar x)$ and 
$\bar{\mathcal F}^{(0)}(x,\bar x)$ are, respectively, 
a $C^\infty$ section of ${\mathcal L}^2$ and a $C^\infty$ section of 
$\bar {\mathcal L}^2$ \cite{St1}. 
The extension of ${\mathcal F}^{(0)}(x,\bar x)$ to genus one was introduced 
in \cite{BCOV1} by the $t$-$t^*$ equation,
\[
\pd_x \pd_{\bar x}{\mathcal F}^{(1)}(x,\bar x) = 
\frac{1}{2}C_{xxx}C_{\bar x\bar x\bar x}e^{2K}G^{x\bar x}G^{x\bar x}-
(\frac{\chi}{24}-1)G_{x\bar x} \;\;.
\]
Geometrically ${\mathcal F}^{(1)}(x,\bar x)$ is understood to represent 
a certain Hermitian norm ( `Quillen's norm' or analytic torsion) of 
a holomorphic line bundle \cite{BCOV1} ( see also \cite{BT, FZY}) over 
the complex structure moduli space. 
The higher genus generalization ${\mathcal F}^{(g)}(x,\bar x) \; (g \geq 2)$ 
are defined by a kind of recursion relation,  
the BCOV holomorphic anomaly equation,
\begin{equation}
\pd_{\bar x} {\mathcal F}^{(g)}=\frac{1}{2}C_{\bar x\bar x\bar x}e^{2K}G^{x\bar x}
G^{x\bar x} 
\Big\{ D_x D_x {\mathcal F}^{(g-1)}+\sum_{r=1}^{g-1} D_{x} {\mathcal F}^{(g-r)} 
D_{x} {\mathcal F}^{(r)} \Big\} \;,
\mylabel{eqn:BCOVanomaly}
\end{equation}
for $C^\infty$ sections ${\mathcal F}^{(g)}(x,\bar x)$ of 
${\mathcal L}^{2-2g}$, more precisely `$C^\infty$-$\Gamma$-modular forms' 
of weight $(2-2g,0)$. 

Recent progresses made in \cite{ABK, GNP} clarify the meaning 
of the anomaly equation (\ref{eqn:BCOVanomaly}) using the wave function 
interpretation of the topological string amplitude \cite{OSV}.
In particular, in \cite{ABK}, similarities of the `$C^\infty$-$\Gamma$-modular 
forms' to the quasi-modular forms in elliptic curves \cite{KZ} has been 
made explicit.

\vskip0.5cm

The general solutions of the BCOV anomaly equation have been obtained by 
certain Feynman rules in \cite{BCOV2}. To present the result, let us 
introduce the notation
$
{\mathcal F}^{(g)}_{r}=\underbrace{D_x \cdots D_x}_{r} {\mathcal F}^{(g)} 
$ 
and define ${\mathcal F}^{(g)}_{r;s}$ 
recursively by 
\[
{\mathcal F}^{(g)}_{r;s+1}=(2g-2+r+s){\mathcal F}^{(g)}_{r;s} \;\;
\;\; ({\mathcal F}^{(g)}_{r;0}={\mathcal F}^{(g)}_r) \;\;,
\]
with the conditions,
\[
{\mathcal F}^{(0)}_{r;1}=0 \;\; (r \leq 2) \;\;;\;\;
{\mathcal F}^{(1)}_{0;1}=\frac{\chi}{24}-1 \;\;,\;\; 
{\mathcal F}^{(1)}_{0;0}=0 \;\;. 
\]
Define perturbative interaction function $P(J,\phi)$ and the source 
function $G(J,\phi)$ by 
\[
P(J,\phi)=\sum_{g \geq 0} \sum_{r,s \geq 0} \lambda^{2g-2} 
{\mathcal F}^{(g)}_{r;s} \frac{J^r}{r!} \frac{\phi^s}{s!}\;\;,\;\;
G(J,\phi)=e^{
-\lambda^2( \frac{1}{2}\,S^{xx} J^2- S^x J \phi -\frac{1}{2}\,S\phi^2 )},
\]
where $\lambda$ is a parameter (string coupling constant) and 
$S^{xx},S^x,S$ represent the propagators determined by integrating 
$e^{2K}D_{\bar x} D^x D^x \bar{\mathcal F}^{(0)}=\pd_{\bar x} S^{xx}$ 
and similar relations for $S^x$ and $S$, see Appendix (A-2). 
One may solve these propagators in the following form, 
\begin{equation}
\begin{aligned}
S^{xx}&=\frac{1}{C_{xxx}}(2\,K_x-\Gamma_{xx}^{\;\; x}+\frac{1}{v^x}\pd_x v^x )
\,,\;
S^x=\frac{1}{2}D_x S^{xx}+\frac{1}{2}(S^{xx})^2\,C_{xxx}+ H^x_1, \\
S &=H^x_1 K_x +\frac{1}{2}D_x S^x + \frac{1}{2}S^{xx}S^{x}C_{xxx}+H_2,
\end{aligned}
\mylabel{eqn:propagators-S}
\end{equation}
where $v^x(x), H_1^x(x)$ represent some (rational) vector fields and $
H_2(x)$ is a rational function on the moduli space. 
These propagators are $C^\infty$ sections of ${\mathcal L}^{-2}$ with 
suitable tensor indices. 
We introduce the holomorphic (meromorphic) functions $f_g(x)$ on the moduli 
space to represent the `constants' of the integration of the anomaly equation 
(\ref{eqn:BCOVanomaly}). Then the solutions of the anomaly equation can be  
formulated in the following perturbative expansion;
\[
e^{-\sum_{g}\lambda^{2g-2}f_g} = 
e^{P(\frac{\pd \;}{\pd J},\frac{\pd \;}{\pd \phi})} G(J,\phi) \Big|_{J=\phi=0} \;.
\]
The logarithm of the right hand side represents summing over 
connected Feynman diagrams with the interaction terms determined by 
$P(J,\phi)$, and we see the perturbative expansion of ${\mathcal F}^{(g)}$ 
at the coefficient of $\lambda^{2g-2}$(, see (6.16) in \cite{BCOV2}). 

\vskip0.3cm

For convenience, we write the resulting expression at the coefficient 
$\lambda^{2g-2}=\lambda^2$,
\[
\begin{aligned}
{\mathcal F}^{(2)}=& 
\frac{5}{24} \, (S^{xx})^3 \, ({\mathcal F}^{(0)}_{3})^2 
- \frac{1}{8} \, (S^{xx})^2 \, {\mathcal F}^{(0)}_4 
-\frac{1}{2}\, (S^{xx})^2 \, {\mathcal F}^{(0)}_{3} {\mathcal F}^{(1)}_{1} 
+\frac{1}{2}\, S^{xx} \, ({\mathcal F}^{(1)}_{1})^2 \\
&
+ \frac{1}{2} \, S^{xx} \, {\mathcal F}^{(1)}_{2} 
+ \frac{\chi}{24} \, S^x \, {\mathcal F}^{(1)}_1 
-\frac{\chi}{48} \, S^x S^{xx} {\mathcal F}^{(0)}_3 
+\frac{\chi}{24}\,(\frac{\chi}{24}-1)\, S + f_2 \;,
\end{aligned}
\]
where by definition ${\mathcal F}^{(0)}_{3}=C_{xxx}$, and $f_2=f_2(x)$ is the 
holomorphic ambiguity. In general, ${\mathcal F}^{(g)}$ is an element of 
$\Gamma_{\infty}({\mathcal L}^{2-2g})$ and may be expressed by 
\begin{equation}
{\mathcal F}^{(g)}(x,\bar x)=
\Gamma(S^{xx},S^x,S;{\mathcal F}_r^{(h<g)}(x,\bar x))+f_g(x)\;\;,
\mylabel{eqn:sol-BCOV}
\end{equation}
where $\Gamma$ represents symbolically the summation over the Feynman 
diagrams.

\vskip1cm \noindent
{\bf (3-3) ${\tt F}_g(t)$ from the topological limit.} 
Following \cite{BCOV1}, we define the `topological limit' of 
(\ref{eqn:sol-BCOV}).
First, the data of the topological limit consists of the normalized solutions 
$w_0(x)$ and $w_1(x)$ at the degeneration point, which determines 
the mirror map $t=t(x)$, and also the initial data for $g=0,1$,
\[
F^{(0)}_3(x)=C_{xxx} \;\;,\;\; 
F^{(1)}_1(x)=\pd_x F^{(1)}(x)\;\;,
\]
where $C_{xxx}$ is the Yukawa coupling (\ref{eqn:Cxxx}) and $F^{(1)}(x)$ 
is the BCOV formula (\ref{eqn:F1}). Then the topological `limit' 
is defined by the following replacements, 
\begin{equation}
G_{x\bar x} \rightarrow \frac{d t}{d x} \frac{d \bar t}{d \bar x} 
\;\;,\;\;
K_x \rightarrow - \pd_x \log w_0(x) \;\;,\;\;
{\mathcal F}^{(g)}(x,\bar x) \rightarrow F^{(g)}(x) \;\;,
\mylabel{eqn:top-limit}
\end{equation}
in the solution (\ref{eqn:sol-BCOV}), which gives 
\begin{equation}
F^{(g)}(x)=\Gamma(S^{xx}(x),S^x(x),S;F_r^{(h<g)}(x))+f_g(x)\;\;.
\mylabel{eqn:top-recursion}
\end{equation}
This is a recursion relation that determines the holomorphic prepotentials 
$F^{(g)}(x)$ as the holomorphic sections of ${\mathcal L}^{2-2g}$ 
starting with the initial data $F^{(0)}_3(x)$ and $F^{(1)}_1(x)$ above. 
Leaving aside 
the holomorphic ambiguity $f_g(x)$, the holomorphic prepotential gives 
the Gromov-Witten potential by
\begin{equation}
\begin{aligned}
{\tt F}_g(t)&=(w_0(x))^{2g-2}\,F^{(g)}(x)  \\
&=
(w_0(x))^{2g-2}\,\Gamma(S^{xx},S^x,S;F_r^{(h<g)}(x))+ 
(w_0(x))^{2g-2}\,f_g(x) \;\;.
\end{aligned}
\mylabel{eqn:GW-Fg}
\end{equation}
The meaning of the topological limit has been discussed recently 
\cite{ABK, GNP} in terms of the 
wave function interpretation of $exp( \sum_{g \geq 0} \lambda^{2g-2} 
{\mathcal F}^{(g)} )$ in \cite{OSV}, however the connection of the 
holomorphic potential $F^{(g)}(x)$ to the Gromov-Witten potential 
${\tt F}_g(t)$ above is still opened mathematically 
(cf. the so-called `mirror theorem' by \cite{Gi, LLY} for $g=0$).  
To determine the ambiguity $f_g(x)$, we have to invoke some regularity 
arguments for ${\tt F}_g(t)$. This 
restricts the possible form of $f_g(x)$. 
Although the regularity arguments put rather strong 
constraints on the possible forms of the ambiguities, we need some 
`boundary' conditions to fix them completely.  We will describe in 
subsection (3-6) the gap conditions at the conifolds which 
has recently introduced in  \cite{HKQ}.

\vskip1cm \noindent
{\bf (3-4) Solving BCOV equation recursively.}  
The general form (\ref{eqn:sol-BCOV}) 
or its topological limit (\ref{eqn:top-recursion}) 
is not so useful for higher genus calculations, since it contains 
the contributions from the large number of connected Feynman diagrams, 
even for $g=4$ or $g=5$. 
On this respect, Yamaguchi and 
Yau\cite{YY} found a nice way to improve the situation. Their idea is to 
formulate a recursion relation for the sections 
$\{\mathcal{ F}^{(g)}(x,\bar{x}) \}$ in the 
form of a differential equation. This avoids the large summation over the 
Feynman diagrams. 

\vskip0.7cm

\noindent
{\bf (3-4.\,a)} 
Following  Yamaguchi and Yau\cite{YY}, let us  
introduce the following expressions,
\begin{equation}
A_{k}=G^{x\bar{x}}\, \theta_x^k \, G_{x\bar{x}} 
\;\;,\;\;
B_k=e^{K(x,\bar x)} \, \theta_x^k \, e^{-K(x,\bar x)} \;\;(k=1,2,\cdots),
\mylabel{eqn:DkA-DkB}
\end{equation}
where $\theta_x=x\frac{d \;}{d x}$. By definition, these satisfy 
\[
\theta_x A_k = A_{k+1}-A_1 A_{k} \;\;,\;\;
\theta_x B_k = B_{k+1}-B_1 B_{k} \;\;.\;\;
\]
Also, since $e^{-K(x,\bar x)}=(\Omega(x),\overline{\Omega(x)})$ 
satisfies the (holomorphic) Picard-Fuchs equation (\ref{eqn:GrPF}) of 
the fourth order, there is a linear relation
\[
B_4 + r_1(x) \, B_3 + r_2(x) \, B_2 + r_3(x) \, B_1 + r_4(x) =0 \;,
\]
with the rational functions $r_k(x)$ which follow from (\ref{eqn:GrPF}). 
Similarly for $A_2(x)$, but from a non-trivial reasoning, we have \cite{YY}
\begin{equation}
A_2=-4\,B_2-2\,B_1(A_1-B_1-1)+\theta_x\log(x C_{xxx}) \; (A_1+2\,B_1+4) + 
r(x) \;\;,
\mylabel{eqn:A2}
\end{equation}
with a rational function $r(x)$, see Appendix (A-3). 
These relations (\ref{eqn:DkA-DkB}) and (\ref{eqn:A2}) entail 
an important property,
\begin{equation}
\theta_x: {\bf C}(x)[A_1,B_1,B_2,B_3] \rightarrow {\bf C}(x)[A_1,B_1,B_2,B_3] \;\;,
\mylabel{eqn:YY-theta}
\end{equation}
i.e., $\theta_x$ acts on the polynomial ring of $A_1,B_1,B_2,B_3$ with the 
coefficients over the rational functions ${\bf C}(x)$. 

\vskip0.7cm

\noindent
{\bf (3-4.\,b)} As for the $\pd_{\bar x}$ operation, it is easy to see  
\[
\pd_{\bar x} : {\bf C}(x)[A_1,B_1,B_2,B_3] 
\rightarrow {\bf C}(x)[A_1,B_1,B_2,B_3][\pd_{\bar x}A_1, \pd_{\bar x} B_1] \;\;.
\]
To show this property, let us note the relations 
$B_2=\theta_x B_1+B_1^2$ and $\pd_{\bar x} B_1 = -x\, G_{x\bar x}$. Then 
for $\pd_{\bar x} B_2$, we have 
\[
\pd_{\bar x} B_2 = -\theta_x (x G_{x\bar x})+2\,B_1 \pd_{\bar x} B_1 
=(1+A_1+2 \, B_1)\,\pd_{\bar x}B_1 \;,
\]
where we use $\theta_x G_{x\bar x}=A_1 \,G_{x\bar x}=
-\frac{1}{x} A_1 \pd_{\bar x} B_1$. Applying $\theta_x$ to this result 
and using $B_3=\theta_x B_2 + B_1 B_2$, we have 
\[
\pd_{\bar x}B_3=(A_2+2\,A_1+3\,B_1+3\,B_2+3\,A_1 B_1+1)\,\pd_{\bar x}B_1 \;.
\]
This shows the claim above.

\newpage
\noindent
{\bf (3-4.\,c)} Now let us focus on the recursion relation  
(\ref{eqn:top-recursion}) for $\mathcal{F}^{(g)}(x,\bar{x})$ 
with the results obtained 
in (3-4.a) and (3-4.b) above. First, we note that the initial conditions are 
given in the ring ${\bf C}(x)[A_1,B_1,B_2,B_3]$ since 
$\mathcal{F}^{(0)}_3(x)=C_{xxx}$ and we have, 
from (\ref{eqn:mathcal-F1}),
\[
\mathcal{F}^{(1)}_1(x,\bar{x})=
\frac{1}{2\, x}\Big\{ -A_1 -(3+h^{11}-\frac{\chi}{12})B_1 
-1-\frac{c_2. H}{12} + \frac{x(57+578\,x-3\,x^2)}{6 \, dis(x) } \Big\} \;.
\]
Also for the propagators we see that $S^{xx}, S^x, S$ belong to 
the ring ${\bf C}(x)[A_1,B_1,B_2,B_3]$, see (\ref{eqn:propagators-S}).
For example, we have 
\begin{equation}
S^{xx}=-\frac{1}{x\,C_{xxx}}( A_1+2\,B_1+4) \;\;,\;\;
S^{x}=\frac{1}{x^2\,C_{xxx}}( 3\,B_1+B_2+2) \;\;.
\mylabel{eqn:S-by-AB}
\end{equation}
The recursion relation 
(\ref{eqn:sol-BCOV}) contains the covariant derivatives $D_x$ 
to define $\mathcal{F}^{(h<g)}_r$ $=
D_x \cdots D_x \mathcal{F}^{(h<g)}(x,\bar{x})$. Note that these  
covariant derivations act inside the ring due to the property 
(\ref{eqn:YY-theta}). Therefore, by induction, we may conclude that 
the prepotentials $\mathcal{F}^{(g)}(x,\bar{x})$ are in the ring 
${\bf C}(x)[A_1,B_1,B_2,B_3]$ for all $g \geq 2$. This 
is the polynomiality found in \cite{YY}. 

\vskip0.3cm

Now we proceed to combine the polynomiality with the integration of 
the BCOV anomaly equation (\ref{eqn:BCOVanomaly}). Following \cite{YY}, 
let us introduce 
$P^{(g)}_n \in {\bf C}(x)[A_1,B_1,B_2,B_3]$ ($P_0^{(g)}=P^{(g)}$) by 
\begin{equation}
P^{(g)}_n =(x^3\, C_{xxx})^{g-1}\,x^n D_x^n \mathcal{F}^{(g)}  
\;\;\;(n=0,1,2,\cdots).
\end{equation}
Then it is straightforward to rewrite the BCOV equation as
\[
\pd_{\bar x} P^{(g)}=\frac{1}{2}\pd_{\bar x} (x\, C_{xxx} S^{xx}) \Big\{ 
P_{2}^{(g-1)} + \sum_{r=1}^{g-1} P_1^{(g-r)}\,P_1^{(r)} \Big\} \;.
\]
Both sides of this equation are linear in $\pd_{\bar x}A_1, \pd_{\bar x}B_1$, 
and if we {\it assume} these two are linearly independent, then we have 
\begin{equation}
\begin{aligned}
&2\; \frac{\pd P^{(g)}}{\pd A_1}  -  
\Big( \frac{\pd P^{(g)}}{\pd B_1} + 
\frac{\pd_{\bar x} B_2}{\pd_{\bar x} B_1}\frac{\pd P^{(g)}}{\pd B_2} +  
\frac{\pd_{\bar x} B_3}{\pd_{\bar x} B_1}\frac{\pd P^{(g)}}{\pd B_3} \Big) 
= 0 \;\;, \\
&\frac{\pd P^{(g)}}{\pd A_1}  
=-\frac{1}{2}\,\Big\{ 
P_{2}^{(g-1)} + \sum_{r=1}^{g-1} P_1^{(g-r)}\,P_1^{(r)} \Big\}  \;\;.
\end{aligned}
\mylabel{eqn:diff-Pg-AB}
\end{equation}
The first equation implies that $P^{(g)}$ is a polynomial of essentially 
three variables. This suppresses the length of the polynomial $P^{(g)}$ 
when $g$ becomes large. A nice choice of 
variables that respects the first equation is given in \cite{YY} by 
\begin{equation}
\begin{aligned}
&B_1=u \;\;,\;\;
A_1=v_1-2\,u-1 \;\;,\;\;
B_2=v_2+u\,v_1 \;\;, \\
&
B_3=v_3+u \Big( 2\,v_1 + 
    \theta_x \log(x C_{xxx})\,v_1-v_2+3\theta_x \log(x C_{xxx}) 
    + r(x) -1 \Big) \;\;.
\end{aligned}
\mylabel{eqn:AB-uv}
\end{equation}
Note that the inverse relation to this may be found easily because the 
above relation is of `upper triangular form'. Using 
the new variables for the first equation of (\ref{eqn:diff-Pg-AB}), 
we have $\frac{\pd \;}{\pd u}P^{(g)}=0$ and conclude, 
\[
P^{(g)} \in {\bf C}(x)[v_1,v_2,v_3]  \; \subset \; 
{\bf C}(x)[u, v_1,v_2,v_3]={\bf C}(x)[A_1,B_1,B_2,B_3] \;\;.
\]
Furthermore, note that the both sides of the second equation 
in (\ref{eqn:diff-Pg-AB}) are polynomial in $u$ of degree less than three. 
Then, writing  
$\frac{1}{2} \{ 
P_{2}^{(g-1)} + \sum_{r=1}^{g-1} P_1^{(g-r)}\,P_1^{(r)} \} 
=: Q_0+u\,Q_1+u^2\,Q_2$, we have
\begin{equation}
\frac{\pd P^{(g)}}{\pd v_1}=-Q_0 \;\;,\;\;
\frac{\pd P^{(g)}}{\pd v_2}=Q_1+(2+\theta_x\log(x C_{xxx}))\,Q_2 \;\;,\;\;
\frac{\pd P^{(g)}}{\pd v_3}= Q_2 \;\;.\;\;
\mylabel{eqn:diff-Pg}
\end{equation}
This is the equation we can solve recursively with the initial 
data $P^{(0)}_3=1$ and $P^{(1)}_1$.

\vskip0.7cm
\noindent
{\bf (3-4.\,d)} The holomorphic ambiguity $f_g$ in (\ref{eqn:sol-BCOV}) 
corresponds to the `constants' of the integration of the differential 
equation (\ref{eqn:diff-Pg}). To make the correspondence more precise, we note 
that $f_g$ in (\ref{eqn:sol-BCOV}) may be identified by the vanishing 
limit of the propagators, i.e. $\mathcal{F}^{(g)} \rightarrow f_g$ when 
$S^{xx},S^x,S \rightarrow 0$. 
Now assume that $P^{(g)} \in {\bf C}(x)[v_1,v_2,v_3]$ is a solution of the 
differential equation (\ref{eqn:diff-Pg}). We substitute in 
$(x^3 C_{xxx})^{1-g}\,P^{(g)}(v_1,v_2,v_3)$ the expressions 
for $v_1,v_2,v_3$ in terms of the propagators, which follow from 
(\ref{eqn:S-by-AB}) and (\ref{eqn:AB-uv}).  
Then the vanishing limit of the propagators gives  
the holomorphic ambiguity $f_g$. In other words, we may write 
\begin{equation}
\mathcal{F}^{(g)}=(x^3 C_{xxx})^{1-g} P^{(g)} + f_g(x) \;\;,
\mylabel{eqn:Fg-by-Pg}
\end{equation}
where we fix the integration `constant' in $P^{(g)}$ by 
the property $P^{(g)}(v_1,v_2,v_3) \rightarrow 0$ when $S^{xx},S^{x},S
\rightarrow 0$.

\vskip1cm \noindent
{\bf (3-5) Relating the topological limits.} 
Let us note that the topological limit (\ref{eqn:top-limit}) with the 
data $w_0(x), w_1(x), t=t(x)$ corresponds to the replacements
\[
A_{1}\to\Big(\frac{d x}{d t} \Big) \theta_x \Big(\frac{d t}{d x}\Big) 
\;\;,\;\;
B_k\to\frac{1}{w_0(x)}\theta_x^k w_0(x) \;\;(k=1,2,3)\;,
\]
in the polynomial solutions ${\mathcal F}^{(g)}=
{\mathcal F}^{(g)}(A_1(x,\bar x),B_k(x,\bar x),x)$. 
We denote the resulting holomorphic potential $F^{(g)}(x)$. 

Now we define $\tilde {\mathcal F}^{(g)}(z,\bar z)$ to be the solutions 
of the BCOV equation in $z$-coordinate with the initial conditions 
$\tilde {\mathcal F}^{(1)}_1(z, \bar z)$ and $\tilde 
{\mathcal F}^{(0)}_3=D_{z}D_{z}D_{z}\tilde {\mathcal F}^{(0)}(z,\bar z)$.  
Since the initial data, in particular for $g=0$, are related by 
\[
\tilde{\mathcal F}^{(0)}_3(z,\bar z)=C_{zzz}(z)=C_{xxx}(\tfrac{1}{z})
\Big( \frac{d x}{d z} \Big)^3 =
{\mathcal F}^{(0)}_3(\tfrac{1}{z},\tfrac{1}{\bar z}) 
\Big( \frac{d x}{d z} \Big)^3 \;\;,
\]
we see that $\tilde {\mathcal F}^{(g)}(z,\bar z)$ and 
${\mathcal F}^{(g)}(x,\bar x)$ are in the same coordinate patch of a  
trivialization of the line bundle ${\mathcal L}$. Hence we have
\begin{equation}
\tilde{\mathcal F}^{(g)}(z,\bar z)= 
{\mathcal F}^{(g)}(\tfrac{1}{z}, \tfrac{1}{\bar z}) \;\;,
\mylabel{eqn:Fz-Fx}
\end{equation}
for the $C^{\infty}$ sections of ${\mathcal L}^{2-2g}$. 
Then, by the data $\tilde{w}_0(z), \tilde{w}_1(z), \tilde{t}=\tilde{t}(z)$ 
given in (\ref{eqn:w0(z)}), the topological limit of  
$\tilde {\mathcal F}^{(g)}(z,\bar{z})=
{\mathcal F}^{(g)}(A_1(\tfrac{1}{z},\tfrac{1}{\bar z}),
B_k(\tfrac{1}{z},\tfrac{1}{\bar z}),\tfrac{1}{z})$ may be achieved by
\[
\begin{aligned}
&A_1(\tfrac{1}{z},\tfrac{1}{\bar z})=
\Big(\frac{dx}{dz}\frac{d\bar x}{d\bar z} G^{z\bar z} \Big)
(-\theta_z) \Big( \frac{dz}{dx} \frac{d \bar z}{d\bar x} G_{z\bar z} \Big) 
\rightarrow 
-\Big(\frac{d z}{d \tilde{t}} \Big) 
\theta_z \Big(\frac{d \tilde{t}}{d z}\Big)-2 \, ,
\\
&B_k(\tfrac{1}{z},\tfrac{1}{\bar z})=
e^{\tilde K(z,\bar z)}(-\theta_z)^k e^{-\tilde K(z,\bar z)} 
\rightarrow
\frac{1}{\tilde{w}_0(z)}(-\theta_z)^k \tilde{w}_0(z) \;, \;\;
(k=1,2,3)~,
\end{aligned}
\]
where the relations $G_{x\bar x}(\frac{1}{z},\frac{1}{\bar z})=
\frac{dz}{dx} \frac{d \bar z}{d\bar x} G_{z\bar z}(z,\bar z)$, 
$K(\frac{1}{z},\frac{1}{\bar z})=\tilde K(z,\bar z)$ have been used. 
We denote the resulting holomorphic potential $\tilde F^{(g)}(z)$. 

According to \cite{BCOV2}, we finally obtain the 
Gromov-Witten potentials 
for $X$ and $X'$ by
\begin{equation}
{\tt F}_g(t)=
(w_0(x))^{2g-2}\,F^{(g)}(x) \;\;,\;\;
\tilde {\tt F}_g(\tilde t)=(\tilde w_0(z))^{2g-2}\, 
\tilde F^{(g)}(z) \;\;,
\mylabel{eqn:GWFxz}
\end{equation}
with the mirror maps $t=t(x)$ and $\tilde t=\tilde t(z)$, respectively.

\vskip0.5cm 
We remark that if we require ${\tt F}_g(t)$ and $\tilde {\tt F}_g(\tilde t)$ 
are regular at $x=0$ and $z=0$, respectively, then the relation 
(\ref{eqn:Fz-Fx}) restricts possible behaviors of the holomorphic (rational) 
function $f_g(x)$, near $x=0$ and $\infty$. Taking these reguarlity 
constraints into accounts, following \cite{BCOV2}, 
we may set the following anzatz for $f_g$,
\begin{equation}
\begin{aligned}
f_g(x)&=a_0+a_1x+\cdots+a_{2g-2}x^{2g-2} \\
&+\frac{b_0+b_1 x + \cdots +b_{2g-3} x^{2g-3}}{(x-3)^{2g-2}} + 
\frac{c_0+c_1 x + \cdots +c_{6g-7} x^{6g-7}}{dis(x)^{2g-2}} \;,
\end{aligned}
\mylabel{eqn:fg-form}
\end{equation}
where $dis(x)=1-57\,x-289\,x^2+x^3$. Alghough $x=3$ does not corresponds 
to any degeneration of the mirror family, we introduce $b_0,\cdots,b_{2g-3}$ 
in this general form (, see section 5 for more detailed analysis on this). 
In this form, we see $10(g-1)+1$ unknown 
parameters which grow linearly in $g$.

\vskip1cm \noindent
{\bf (3-6) The gap conditions at conifolds.} 
One of the most subtle parts in solving the BCOV anomaly equation is 
to fix the holomorphic ambiguities $f_g(x)$ whose general form 
has been argued in (\ref{eqn:fg-form}).  To determine the unknown constants 
contained in $f_g(x)$, 
we may use the first few terms of $N_g(d)$ in the expansion (\ref{eqn:Fg-def}) 
if they are known from other methods, e.g. enumerative geometry. In many 
cases, one may expect $n_g(d)=0$ for lower $d$ assuming that 
$n_g(d)$ counts the number of genus $g$ curves in $X$ of degree $d$ 
and also some genus formula for curves, see e.g. \cite{KKV}. 
However these conditions are not sufficient to determine $f_g(x)$ 
in general, and this fact reduces the predictive power of the BCOV 
equation for determining the Gromov-Witten potentials ${\tt F}_g(t)$. 
Recently, on this problem, Huang, Klemm and Quackenbush \cite{HKQ} 
have found that a certain 
vanishing property (the gap condition) at conifolds provides  
considerably strong conditions for $f_g(x)$. The gap condition has been 
tested for quintic hypersurface in ${\bf P}^4$ and other cases that have 
the mirror family over ${\bf P}^1$ with only one conifold singularity. 

\vskip0.3cm
The gap condition in \cite{HKQ} arises from the topological limit 
made around a conifold singularity. Let $x=c$ be a conifold singularity 
of the mirror family, or the corresponding singularity of the Picard-Fuchs 
differential equation. In our case, $c$ may be one of the three singularities 
$\alpha_1,\alpha_2,\alpha_3$ in (\ref{eqn:P-indices}). As we observe in 
(\ref{eqn:P-indices}), the indices $\rho_k$ at the conifold are all 
integral but have one degeneracy, which indicates there exists one solution 
with logarithmic singularity. 

Assume $(\rho_1,\rho_2,\rho_3, \rho_4)=(0,1,1,2)$, and normalize the 
logarithmic solution \break 
$\log(s) w_1^c(s) + O(s^1) $ 
by requiring $w_1^c(s)= s+O(s^2)$ $(s=(x-c))$. 
Then, according to the Picard-Lefschetz theory, the series $w_1^c(s)$ 
represents the (normalized) period integral of the vanishing cycle. 
$w_1^c(s)$ together with the logarithmic solution 
corresponds to the indices $\rho_2=\rho_3=1$. For the index 
$\rho_4=2$ we have the solution of the form $w_2^c(s)=s^2+ O(s^3)$. Then,  
making a suitable linear combination with $w_1^c(s)$ and $w_2^c(s)$, we 
may fix the solution for the index $\rho_1=0$ by the property 
\[
w_0^c(s)=1+ O(s^3)  \;\;.
\]
By the data of the topological limit at the conifold $x=c$, we mean 
the series data $w_0^c(s), w_1^c(s)$ with the `mirror map' 
$s=s(U)$ defined by 
\[
k_U U=\frac{w_1^c(s)}{w_0^c(s)} \;\;,
\]
where $k_U$ is a constant characterized below. 

\vskip0.3cm

The gap condition arises from the topological limit 
${\mathcal F}^{(g)}_c(s,\bar s) \rightarrow F^{(g)}_c(s)$ at each conifold. 
We define this topological limit, 
in the exactly same way as described in subsections (3-3),(3-5), 
by the replacements 
\[
A_1(s+c,\bar s+\bar c) \rightarrow (s+c) \frac{d \; }{d s} \log\frac{d U}{ d s} \;\;,\;\;
B_k \rightarrow 
\frac{1}{w_0^c(s)} \Big( (s+c)\frac{d \;}{d s} \Big)^k 
w_0^c(s) \;\;.
\]
in the relation ${\mathcal F}_c^{(g)}(s,\bar s)=
{\mathcal F}^{(g)}(A_1(x,\bar x),B_k(x,\bar x),x)$.

The observation made in \cite{HKQ} based on the physical interpretation 
of the vanishing cycles \cite{St2} is the following: {\it There exists a choice 
of the constant $k_U$, under which we have
\begin{equation}
{\tt F}_c^{(g)}(U)=(w_0(s))^{2g-2}\, F_c^{(g)}(s)=
\frac{|B_{2g}|}{2\,g\,(2\,g-2)} \frac{1}{U^{2g-2}} + O(U^0) \;,
\mylabel{eqn:gap-condition}
\end{equation}
for $g\geq 2$ (and ${\tt F}_c^{(1)}(s)=-\frac{1}{12} \log U + O(U^0)$). }
Since the leading behavior $F_c^{(g)}(s) \sim 
\frac{\text{const.}}{U^{2g-2}} + \cdots $ 
can be verified in general, the above equation provides $(2g-2)-1$ vanishing 
conditions for the coefficients of $\frac{1}{U^k} \;(1\leq k \leq 2g-3)$, 
the {\it gap condition}. Note that once we find $k_U$ at some $g$, then 
the leading term in (\ref{eqn:gap-condition}) provides an additional 
condition for each other value of $g$.
It has been observed for the quintic and similar 
Calabi-Yau threefolds \cite{HKQ} that these vanishing conditions provides 
an efficient way to determine the holomorphic 
ambiguity $f_g(x)$ for higher values of $g$.

\newpage
\section{{\bf Calculations}}

Here we present some details of our calculations of the Gromov-Witten 
potentials ${\tt F}^X_g(t)={\tt F}_g(t)$ and 
${\tt F}^{X'}_g(\tilde t)=\tilde {\tt F}_g(\tilde t)$, and list the 
resultant Gopakumar-Vafa invariants, $n^X_g(d)$ and $n^{X'}_g(d)$ for 
$g \leq 5$ in Tables 1 and 2. 

\vskip0.5cm
\noindent
{\bf (4-1) Expansions about the conifolds.} 
The evaluations of the Gromov-Witten potentials 
\[
{\tt F}_g(t)=(w_0(x))^{2g-2} F^{(g)}(x) \;\;,\;\;
\tilde {\tt F}_g(\tilde t)=(\tilde w_0(z))^{2g-2} \tilde F^{(g)}(z) \;\;,\;\;
\]
are straightforward with the topological data $w_0(x),w_1(x), t=t(x)$ and 
$\tilde w_0(z),\tilde w_1(z)$, $\tilde t=\tilde t(z)$ as described precisely 
in the previous sections. 
For the expansion about the conifolds, however, 
we need to make the series expansions about 
$x=\alpha_k$ $(k=1,2,3)$ given by the algebraic equation 
$1-57\,x-289\,x^2+x^3=0$. To achieve this, 
we first write the Picard-Fuchs equation 
\begin{equation}
\sum_{k=0}^4 p_k (\alpha,s) \Big( \frac{d \;}{d s} \Big)^k w^c(s) =0 
\mylabel{eqn:PF-a}
\end{equation}
in the coordinate $s=x-\alpha$ with some polynomials $p_k(\alpha,s)$. 
Note that $\alpha$ may be taken to be any of $\alpha_k$ since we only need  
the relation $1-57\,\alpha-289\,\alpha^2+\alpha^3=0$ in the derivation. 
Now we try to find the solutions of the form 
\[
w^c(\alpha,s)=\sum_{n \geq 0} c_n(\alpha) \, s^{n+\rho} \;\;,
\]
for each choice of the index $\rho=0,1,1,2$. Namely we solve the 
differential equation over the ring 
${\mathcal R}_\alpha={\bf C}[\alpha]/(\alpha^3-289\,\alpha^2-57\,\alpha+1)$. 
Solving the Picard-Fuchs equation (\ref{eqn:PF-a}) over ${\mathcal R}_\alpha$ 
is rather technical, but turns out quite useful since we can impose the gap 
conditions at the three conifold points $\alpha_k$ at one time. 

Recall that the gap conditions may be imposed by making the data 
$w_0^c(\alpha,s)$, $w_1^c(\alpha,s)$ and $s=s(U)$ as defined in the 
subsection (3-6). After some calculations, for the solutions, we obtain  
\[
\begin{aligned}
w_0^c(\alpha,s)&=1
- \Big( 
{\frac {82833753}{33614}} 
+{\frac {1555547739}{134456}}\,\alpha
-{\frac {16148435}{403368}}\,{\alpha}^{2}
\Big) {s}^{3} + \cdots  \\
w_1^c(\alpha,s)&=s - 
\Big( 
 {\frac {64163}{1372}} 
+{\frac {83161}{343}}\,\alpha
-{\frac{1151}{1372}}\,{\alpha}^{2} 
 \Big) {s}^{2} 
+ \cdots  \;\;,
\end{aligned}
\]
and also, inverting the defining relation 
$k_U U=\frac{w_1^c(\alpha,s)}{w_0^c(\alpha,s)}$, we have 
\[
s(U)=
k_U U+ \Big( 
{\frac {64163}{1372}}
+{\frac {83161}{343}}\,\alpha 
-{\frac {1151}{1372}}\,{\alpha}^{2}
\Big) (k_U {U})^{2}+ \cdots  
\]
Using the above data, we can evaluate the holomorphic potential 
$F_c^{(g)}(s)$ in the following form,
\[
{\tt F}_c^{(g)}(U)= 
\frac{R_{2g-2}(\alpha)}{(k_U U)^{2g-2}} + 
\frac{R_{2g-1}(\alpha)}{(k_U U)^{2g-3}} + \cdots +  
\frac{R_{1}(\alpha)}{(k_U U)} + O(U^0) \;,
\]
with $R_{k}(\alpha)=c_{k,2} \, \alpha^2 +c_{k,1} \,\alpha + c_{k,0}$. 
Since $1, \alpha, \alpha^2$ are linearly independent, 
the gap condition (\ref{eqn:gap-condition}) entails $3  (2g-3)$ 
conditions, or $3 (2g-2)$ conditions once $k_U$ is fixed.  Thus we can 
impose the gap conditions at the three conifold points at once 
in this algebraic manipulation. 

\vskip1cm
\noindent
{\bf (4-2) Examples ($g=2,3$).}  
We use the gap condition above extensively together 
with some natural vanishing assumptions to fix the $10\,g-9$ unknown 
parameters in $f_g(x)$, see (\ref{eqn:fg-form}). Here we illustrate how 
we impose the additional vanishing conditions using the cases $g=2$ and 
$g=3$. For $g=2$, we have to fix $10 \,g-9=11$ unknown parameters 
among which $3\,(2g-3) = 3$ may be determined from the gap conditions. 
To fix the remaining $8$ parameters, we note the following $g=1$ 
Gopakumar-Vafa invariants which follow from the BCOV formula (\ref{eqn:F1});
\[
\begin{array}{c | c c c c c c}
         & 1 & 2 & 3 & 4 & 5   & 6 \\
\hline
n_1^X(d) & 0 & 0 & 0 & 0 & 588 & \cdots \\
n_1^{X'}(d) & 0 & 0 & 196 & 99960 & 34149668 & \cdots\\
\end{array}
\]
{}From the higher genus calculations done in several examples,  
see \cite{HST1, HKQ} for example, 
we observe that the vanishing $n_{g-1}(d)=0$ indicates $n_g(d)=0$. 
This observation seems to be a natural consequence of the 
geometrical meaning of the Gopakumar-Vafa invariants that $n_h(d)$ 
is evaluating the Euler numbers 
of the degeneration loci in the genus $g$ curve of degree 
$d$ \cite{GV, KKV, HST2}. 
Assuming that this vanishing condition holds in our case, we have 
\[
n_2^X(d)=0 \;\; (d=1,\cdots,4),\;\;
n_2^{X'}(d)=0  \;\; (d=1,2), \;\;
n_2^{X}(0)=n_2^{X'}(0)=\frac{\chi}{5760} \;\;,
\]
which provide $8$ conditions sufficient to fix 
$f_2(x)$. Using these conditions 
we obtain for the holomorphic potential $F^{(2)}(x)$,  
\[
\begin{aligned}
F^{(2)}(x)= & (x^3\,C_{xxx})^{-1}
\Big( 
\frac{2989}{288} \,v_3 + \frac{49}{24}\, v_1\,v_2 - \frac{5}{24} v_1^3 + 
\frac{p_2(x)}{(x-3) dis(x)}\,v_2  \\
& + \frac{p_{1,1}(x)}{(x-3)\, dis(x)} \,v_1^2 
+\frac{p_1(x)}{(x-3)\, dis(x)^2} v_1 + 
\frac{p_3(x)}{(x-3)\,dis(x)^2 }  \Big) + f_2(x) \;\;,
\end{aligned}
\]
with some polynomials $p_1(x),p_2(x), p_3(x), p_{1,1}(x)$, which we 
leave implicit, and 
\[
\begin{aligned}
f_2(x)&=
  - \frac{359293}{2520}    
  + \frac{1850909\,x}{20160} - \frac{2081\,x^2}{6720} 
- \frac{15739}{24\,{(x -3 ) }^2} 
+ \frac{38147}{84\,( x-3 ) } \\
&
+\frac{1}{dis(x)^2}\Big( 
\frac{264137}{720} - \frac{1881913}{45}x + \frac{39189063}{40}x^2 + 
   \frac{72541963}{6}x^3  \\
& \hskip5cm + \frac{7353789043}{240}x^4 - 
   \frac{8892629}{90}x^5 \Big) \;\;.
\end{aligned}
\]
Also the leading term of the conifold expansion $F_c^{(2)}(s) = 
\frac{1}{240} \frac{1}{U^2} + \cdots$ determines the constant $k_U$ by  
\[
k_U^2=240\, \Big( 
\frac{1183163}{1120}\,\alpha^2 + \frac{58293}{280}\,\alpha 
- \frac{4091}{1120} \Big) \;\;.
\]
The resultant Gopakumar-Vafa invariants $n_2^X(d)$ and $n_2^{X'}(d)$ 
are listed in Table 1 and Table 2. 

For $g=3$ calculation, since $k_U$ has been fixed as above, we have 
$3\,(2g-2) = 12$ constraints from the gap condition to fix $10\,g-9=
21$ parameters in $f_3(x)$. Fortunately, we have enough additional 
vanishing conditions from the $g=2$ results; 
$n_2^{X}(d)=0 \;\;(d=1,\cdots,7)\;$, $n_2^{X'}(d)=0 \;\;(d=1,\cdots,4)$, 
see Table 1 and Table 2,
We may adopt the following $9$ conditions 
\[
n_3^X(d)=0 \;\; (d=1,\cdots,5) \;,\;\;
n_3^{X'}(d)=0  \;\;(d=1,2)\; , \;\;
n_3^{X}(0)=n_3^{X'}(0)=\frac{-\chi}{1451520} \;\;,
\]
to fix $f_3(x)$. 

We have continued this process up to $g=5$.  Although we may continue 
this further to higher $g$, the exact value of $g$ where this process 
might break down is not clear to us (, see the discussion in 
the next section).

\vskip1cm

\section{{\bf Conclusion and discussions}}

We have determined the Gromov-Witten potentials ${\tt F}^X_g$ and 
${\tt F}^{X'}_g$, up to $g=5$, of the Grassmannian and the 
Pfaffian Calabi-Yau threefolds using the mirror symmetry. 
Our calculations are based on the original BCOV holomorphic anomaly equation 
\cite{BCOV1, BCOV2} and the polynomiality in the solutions found 
in \cite{YY}. In particular, following \cite{HKQ}, we used extensively 
the gap conditions at the conifold singularities to determine the 
holomorphic ambiguities $f_g$. 

Apart from these computational aspects of the Gromov-Witten invariants, 
we have also remarked that the (mirror) Picard-Fuchs differential equation 
has a similar property to that appeared in the mirror symmetry 
of a K3 surface of degree 12. For a K3 surface of degree 12, 
the number of the Fourier-Mukai partners is two 
, i.e. $|FM(X)|=2$ \cite{Og, HLOY1}. 
One may expect a similar result for the Grassmannian and 
the Pfaffian Calabi-Yau manifolds, i.e. 
there is no more variety which is derived equivalent to these up to 
isomorphisms. 
Also one may expect that $X'$ appears as a suitable moduli space of 
stable sheaves on $X$, which is the case for the K3 surfaces 
of degree 12.

\vskip0.5cm

Finally we comment on the singularity we see at $x=3$ in 
(\ref{eqn:P-indices}). This point does not corresponds to a singularity 
of the mirror manifold $Y_x$ in (2-1), see \cite{Ro} for more 
details. In fact, we see from the indices at $x=3$, there is no local 
monodromy around this point. However, we can formulate 
additional `gap condition' which may be used to determine 
the holomorphic ambiguity $f_g$. 
Let us fix the local solutions corresponding to $\rho=0,1,3,4$, 
respectively, by the following properties;
\[
\begin{aligned}
&w_0(s)=1-\frac{s^2}{42}+O(s^5) \;\;,\;\;\;\;
w_1(s)=s-\frac{8}{21}\, s^2 + O(s^5) \;\;,\;\; \\
&w_2(s)=s^3-\frac{191}{210}\,s^4 + O(s^5) \;\;,\;\;
w_3(s)=s^4 + O(s^5) \;\;,
\end{aligned}
\]
where $s=x-3$. Then similarly to the conifold points, 
one may define the topological limit with the data $w_0(s), w_1(s)$ and 
the mirror map $U=\frac{w_1(s)}{w_0(s)}$. Then corresponding to the 
gap condition (\ref{eqn:gap-condition}) at the conifolds, we {\it observe} 
that the following vanishing property,  
\[
{\tt F}^{(g)}(U)=(w_0(s))^{2g-2}\, F^{(g)}(s) = 
0\,\frac{1}{U^{2g-2}}+\cdots+0\,\frac{1}{U}+O(U^0) \;\;, 
\]
holds for $g\leq 5$. Note that by the form $f_g$ in (\ref{eqn:fg-form}) 
this expansion can start from $\frac{1}{U^{2g-2}}$ in general. However 
the ${\tt F}^{(g)}(U)$ is regular as above since there does not appear any 
massless state (or vanishing cycle) at $x=3$. We may utilize this property 
to determine the unknown constants in $f_g$. Thus, together with 
the gap conditions at the conifolds, we have 
$8(g-1)$ conditions in total, and hence 
in order to fix $f_g$ completely we need additionally $2\,g-1$ 
vanishing conditions, $n_g^{X}(d)=0,$ $n_g^{X'}(d')=0$ for lower degrees 
$d$ and $d'$. From the results at $g=5$, one may expect that 
the calculations done in section 4 may be continued to 
considerably higher value of $g$, like the case of the quintic \cite{HKQ}.

\vfill

\[
 \begin{array}{|c | l l l}
\hline
d & g=0 & g=1 & g=2 \cr
\hline
1&  196 & 0 & 0 \cr 
2&  1225 & 0 & 0 \cr 
3&  12740 & 0 & 0 \cr 
4&  198058 & 0 & 0 \cr 
5&  3716944 & 588 & 0 \cr 
6&  79823205 & 99960 & 0 \cr 
7&  1877972628 & 8964372 & 0 \cr 
8&  47288943912 & 577298253 & 99960 \cr 
9& 1254186001124 & 31299964612 & 47151720 \cr 
10& 34657942457488 & 1535808070650 & 7906245550 \cr 
11& 990133717028596 & 70785403788680 & 858740761340 \cr 
12& 29075817464070412 & 3129139504135680 & 73056658523632 \cr 
13& 873796023687033916 & 134357808679487260 & 5317135023839604 \cr 
14& 26782042513523921505 & 5648906799029453044 & 347478656042915187 \cr 
15& 834938101511448746224 & 233816422635171601176 & 20996780173465726448 \cr 
16& 26417440686921151630504 & 9563588497688111378163  
        & 1195726471411561809370 \cr 
17& 846787615783681427068332 & 387581693402348794414352  
        & 65017598161994032437484 \cr 
\end{array}
\]

\vskip0.5cm

\[
\begin{array}{|c|l l l}
\hline
d & g=3 & g=4 & g=5 \cr
\hline
1&  0 & 0 & 0 \cr 
2&  0 & 0 & 0 \cr 
3&  0 & 0 & 0 \cr 
4&  0 & 0 & 0 \cr 
5&  0 & 0 & 0 \cr 
6&  0 & 0 & 0 \cr 
7&  0 & 0 & 0 \cr 
8&  0 & 0 & 0 \cr 
9&  -1176 & 0 & 0 \cr 
10& 325409 & 0 & 0 \cr 
11& 956485684 & -25480 & 0 \cr 
12& 301227323110 & 27885116 & 3675 \cr 
13& 52490228133616 & 67509270780 & 73892 \cr 
14& 6617949361316377 & 28917316111159 & 9783073244 \cr 
15& 676939616238018840 & 6764898614128228 & 13255130550228 \cr 
16& 59768711735781062098 & 1117634949252974670 & 6169573531612148 \cr 
17& 4730781899004364783412 & 146451269357268794212 
          & 1690718304511081104 \cr 
18& 344157075745064476608707 & 16239378567823605642392 
          & 332432097873830811843 \cr 
\end{array}
\]
\vskip0.3cm
\noindent 
{\bf Table 1.} 
Gopakumar-Vafa invariants $n_g^{X}(d)$ $(g \leq 5)$ of the Grassmannian 
Calabi-Yau threefold $X=Gr(2,7)_{1^7}$.

\newpage

\[
\begin{array}{|c|lll}
\hline
d & g=0 & g=1 & g=2 \cr
\hline
1&   588 & 0 & 0 \cr 
2&   12103 & 0 & 0 \cr 
3&   583884 & 196 & 0 \cr 
4&   41359136 & 99960 & 0 \cr 
5&   3609394096 & 34149668 & 12740 \cr 
6&   360339083307 & 9220666238 & 25275866 \cr 
7&   39487258327356 & 2163937552736 & 21087112172 \cr 
8&   4633258198646014 & 466455116030169 & 11246111235996 \cr 
9&   572819822939575596 & 95353089205907736 & 4601004859770928 \cr 
10&  73802503401477453288 & 18829753458134112872 
         & 1586777390750641117 \cr 
11&  9831726718738661469404 & 3632247018393524104896 
         & 486768262807329916336 \cr 
12&  1346383795156980043546418 & 689243453496908009355852 
         & 137262882246594110683614 \cr 
\end{array}
\]

\vskip0.5cm

\[
\begin{array}{|c|lll}
\hline
d & g=3 & g=4 & g=5 \cr
\hline
1&   0 & 0 & 0 \cr 
2&   0 & 0 & 0 \cr 
3&   0 & 0 & 0 \cr 
4&   0 & 0 & 0 \cr 
5&   0 & 0 & 0 \cr 
6&   1225 & 0 & 0 \cr 
7&   22409856 & 0 & 0 \cr 
8&   58503447590 & 25371416 & 3675 \cr 
9&   67779027822044 & 216888021056 & 33575388 \cr 
10&  50069281882780727 & 521484626374894 & 1111788286385 \cr 
11&  27893405899311185184 & 660609023799091444 & 5358750700883104 \cr 
12&  12822179880173592308422 & 568693999386204794172 & 
           11048054952421812976 \cr 
13&  5131002509749249793297316 & 
           377653013301230457157640 & 14053721920121779703948 \cr 
\end{array}
\]

\vskip0.3cm
\noindent
{\bf Table 2.} 
Gopakumar-Vafa invariants $n_g^{X'}(d)$ $(g \leq 5)$ of the Pfaffian  
Calabi-Yau threefold $X'$.

\vskip1cm
\appendix
\section{}

\noindent
{\bf (A-1) The Pfaffians of ${\mathcal S}(y)$.}
The $7 \times 7$ skew symmetric matrix ${\mathcal S}(y)$ parametrized 
by $[1,y,0]$ in the subsection (2-2) has the following form,
\[
{\mathcal S}= \left( 
\begin{matrix} 
0 & - u_3 & -y\, u_4 & 0 & 0 & y \, u_0 & u_1 \\
u_3 & 0 & - u_5 & -y \, u_6 & 0 & 0 & y \, u_2 \\
y\, u_4 & u_5 & 0 & -u_0 & -y\,u_1 & 0 & 0 \\
0 & y\, u_6 & u_0 & 0 & -u_2 & - y\, u_3 & 0 \\
0 & 0 & y \, u_1 & u_2 & 0 & -u_4 & -y \, u_5 \\
-y\, u_0 & 0 & 0 & y\,u_3 & u_4 & 0 & -u_6 \\
-u_1 & -y\, u_2 & 0 & 0&  y\,u_5 & u_6 & 0 \\
\end{matrix}
\right) \;\;,
\]
where $[u_0,\cdots,u_6] \in \vP^6$. Then the explicit form of the 
Pfaffians, $p_k({\mathcal S})$ are 
\[
\begin{aligned}
p_0({\mathcal S})&=
y^3\,u_1u_2u_3-y^2\,(u_3u_5^2+u_1u_6^2)-y\,u_0u_2u_4+u_2u_5u_6 \;\;,\\
p_1({\mathcal S})&=
y^3\,u_3u_4u_5-y^2\,(u_5u_0^2+u_3u_1^2)-y\,u_2u_4u_6+u_0u_1u_4 \;\;,\\
p_2({\mathcal S})&=
y^3\,u_0u_5u_6-y^2\,(u_5u_3^2+u_0u_2^2)-y\,u_1u_4u_6+u_2u_3u_6 \;\;,\\
p_3({\mathcal S})&=
y^3\,u_0u_1u_2-y^2\,(u_2u_4^2+u_0u_5^2)-y\,u_1u_3u_6+u_1u_4u_5 \;\;,\\
p_4({\mathcal S})&=
y^3\,u_2u_3u_4-y^2\,(u_2u_0^2+u_4u_6^2)-y\,u_1u_3u_5+u_0u_3u_6 \;\;,\\
p_5({\mathcal S})&=
y^3\,u_4u_5u_6-y^2\,(u_6u_1^2+u_4u_2^2)-y\,u_0u_3u_5+u_1u_2u_5 \;\;,\\
p_6({\mathcal S})&=
y^3\,u_0u_1u_6-y^2\,(u_1u_3^2+u_6u_4^2)-y\,u_0u_2u_5+u_0u_3u_4 \;\;.
\end{aligned}
\]

\vskip0.3cm
\noindent
{\bf (A-2) Propagators $S^{xx}, S^{x}, S$.} These propagators are 
defined in \cite{BCOV2} by integrating 
\[
e^{2K}D_{\bar x}D^xD^x \bar{\mathcal F}^{(0)}=\pd_{\bar x}S^{xx} \;\;,\;\;
G_{\bar x x}S^{xx} = \pd_{\bar x} S^x \;\;,\;\;
G_{\bar x x} S^x = \pd_{\bar x} S \;\;.
\]
Using the special K\"ahler geometry relation (\ref{eqn:sp-geom-rel}), one may 
easily verify that (\ref{eqn:propagators-S}) solves these equations. The 
explicit forms $v^x(x),H_1^x(x),H_2(x)$ are determined following \cite{BCOV2},
\[
v^x(x)=\frac{1}{x^4} \;\;,\;\;
H_1^x(x)=-\frac{1}{2}\frac{1}{x^2\,C_{xxx}}(12-r(x)) \;\;,\;\;
H_2(x)=-\frac{1}{x}H_1^x(x) \;\;,
\]
where $r(x)$ is the rational function in (\ref{eqn:A2}), see also (A-2) below. 
The topological limits of these propagators in the $z$ coordinate 
have similar forms to those found in \cite{BCOV2} for the quintic,
\[
\begin{aligned}
&
S^{zz}=\frac{1}{C_{zzz}}\pd_z 
\log \Big\{ \Big( \frac{f(z)}{\tilde w_0(z)} \Big)^2 
           \frac{d z}{d \tilde t} \Big\}  \;\;,\;\;
S^z=\frac{1}{C_{zzz}}\Big\{ 
\big(\pd_z \log \frac{f(z)}{\tilde w_0(z)} \big)^2 - 
\pd_z^2 \log \frac{f(z)}{\tilde w_0(z)} \Big\}  \;, \\
&
S=\Big\{ 
S^z-\frac{1}{2}D_zS^{zz}-\frac{1}{2}\big( S^{zz} \big)^2 C_{zzz} \Big\} 
\pd_z  \log \frac{f(z)}{\tilde w_0(z)} + 
\frac{1}{2}D_zS^z + \frac{1}{2}S^{zz}S^z C_{zzz} \;\;,
\end{aligned}
\]
where $f(z)=z$. Rather complicated forms of 
$v^x, H_1^x, H_2$ above have been found from the 
latter expressions of $S^{zz},S^z,S$. 

\vskip0.3cm
\noindent
{\bf (A-3) The derivation of $A_2$ in (\ref{eqn:A2}).} 
The relation (\ref{eqn:A2}) follows from the definitions 
\[
\pd_{\bar x} S^{xx} = e^{2K}(G^{x\bar x})^2 C_{\bar x\bar x\bar x} \;\;,\;\;
\pd_{x} C_{\bar x \bar x \bar x} = 0\;\;,
\]
where $C_{\bar x\bar x\bar x}=D_{\bar x}D_{\bar x}D_{\bar x} 
\bar{\mathcal F}^{(0)}(x,\bar x)$ is the anti-holomorphic Yukawa coupling. 
{}From these relations, after some algebra, we have 
\[
\pd_{\bar x}(x C_{xxx} \,\theta_x S^{xx} ) = 
2x \,\{ K_x - \Gamma_{xx}^{\;\;x} \} \pd_{\bar x} (x C_{xxx} S^{xx}) \;\;.
\]
Now from the special geometry relation (\ref{eqn:sp-geom-rel}), we 
have $\pd_{\bar x}(K_x -\Gamma_{xx}^{\;\;x})=-G_{x\bar x}+C_{xxx} 
\pd_{\bar x} S^{xx}$. Using this relation for 
$\pd_{\bar x}( x C_{xxx} S^{xx} )$ in the right hand side, 
we obtain
\begin{equation}
\begin{aligned}
\pd_{\bar x}(x C_{xxx} \, \theta_x S^{xx} ) & = 
2\,x\,(K_x-\Gamma_{xx}^{\;\;x})\Big\{ \pd_{\bar x} 
\big(x(K_x-\Gamma_{xx}^{\;\;x}) \big) + x\, G_{x\bar x} \Big\}  \\
& = \pd_{\bar x}\Big\{ 
(x K_x-\Gamma_{xx}^{\;\;x})^2+(x K_x)^2-2(\theta_x-1)(x K_x) \Big\}. \\
\end{aligned}
\mylabel{eqn:A2-int}
\end{equation}
We may express this relation in terms of $A_1, B_1,B_2$ and $B_3$ as follws,
\[
\begin{aligned}
\pd_{\bar x}(x C_{xxx} \theta_x S^{xx} )& = 
\pd_{\bar x}\Big( 
-A_2+A_1^2-2\,B_2+2\,B_1^2+\theta_x \log(x C_{xxx}) \, (A_1+2\,B_1+4) \Big) \\
&=\pd_{\bar x}( A_1^2+2\,A_1 B_1 + 2\,B_2 -2\,B_1) \;\;,
\end{aligned}
\]
where, for the first line,  we use the expression 
$S^{xx}=\frac{-1}{x C_{xxx}}(A_1+2\,B_1+4)$ in (3-4.c) and  
the relation $\theta_x A_1=A_2-A_1^2$.  
This determines the form $A_2(x)$ up to a holomorphic (rational) function. 
Substituting the series data (\ref{eqn:w0(x)}) under the topological 
limit (\ref{eqn:top-limit}), we finally find 
\[
r(x)=11 - \frac{36}{7\,(x-3)} -  
 \frac{4\,\left( 10 - 331\,x - 751\,x^2 \right) }{7\,dis(x) } \;\;, 
\]
in the relation (\ref{eqn:A2}).

\vskip1cm

\end{document}